%% file: main.tex
\documentclass{article}
\usepackage{amsfonts}
\usepackage{amssymb}
\usepackage{fullpage}
\usepackage[utf8]{inputenc}
\usepackage[numbers]{natbib}
\usepackage{enumitem}
\usepackage{graphicx}                   
\usepackage{siunitx}                    
\usepackage[cmex10,intlimits]{amsmath}  
\usepackage{amsthm}                     
\usepackage{amssymb}                    
\usepackage{amsfonts}                   
\usepackage{bm}                         
\usepackage{booktabs}                   
\usepackage{array}                      
\usepackage{enumitem}                   
\usepackage{authblk}                    

\usepackage{xpatch}

\xpatchbibmacro{name:andothers}{%
  \bibstring{andothers}%
}{%
  \bibstring[\emph]{andothers}%
}{}{}
\newcommand{\etal}{\textit{et al}.}

\input{preamble.tex}
\date{}
\title{Golden-Thompson  via pinching inequality}
\author{Saket Choudhary\\
saketkc@gmail.com}
\affil[1]{Computational Biology and Bioinformatics\\

University of Southern California}

\date{}

\begin{document}

\maketitle

\begin{abstract}
For two hermitian matrices $A$ and $B$, Golden-Thompson inequality \cite{golden1965lower, thompson1965inequality} states that
$$
\mathrm{tr}\left[ \exp{A+B} \right] \leq \mathrm{tr}\left[ \exp{A}\exp{B} \right]
$$

Multiple proofs exist for this inequality \cite{forrester2014golden}. Though, most of them lack an intuitive argument. Sutter \etal \  \cite{sutter2017multivariate} presented an intuitive proof using spectral pinching. We elaborate on their approach here.
\end{abstract}

\section{Introduction}
Given $A,B$ hermitian matrices, Golden-Thompson inequality \cite{golden1965lower, thompson1965inequality} states that:
$$
\mathrm{tr}[\exp{A+B}] \leq \mathrm{tr}[\exp{A}\exp{B}]
$$

It is trivial if $A, B$ commute i.e. $AB = BA$. A quick example would be $A = \begin{pmatrix} 1 & 0 \\ 0 & -1 \end{pmatrix}$ and $B = \mathbb{I}_2$, while it does not hold for $A = \begin{pmatrix} 1 & 1\\ 0 & 2\\ \end{pmatrix}$ and $B = \begin{pmatrix} 1 & 2\\ 2 & 3\\ \end{pmatrix}$ as $AB \neq BA$.
There have been multiple proofs \cite{forrester2014golden}, but they lack an intuitive motivation. Sutter \etal \  \cite{sutter2017multivariate} presented a more intuitive proof using spectral pinching. We elaborate on their approach here. There are no contributions in terms of novelty here and is completely expository. We believe it might be helpful for new learners in the field.

\section{Spectral Pinching}

Consider a square complex matrix $A$ partitioned as a $r \times r$ block matrix: $A = \begin{pmatrix} A_{11} & A_{12} & \cdots & A_{1r} \\
A_{21} & A_{22} & \cdots & A_{2r} \\
\vdots & \ddots  & \ddots & \vdots \\
A_{r1} & A_{r2} & \cdots & A_{rr} \\
\end{pmatrix}$ . 
We can decompose this matrix intro two matrices comprising the diagonal and the off-diagonal elements respectively,
\begin{align*}
    A &= A_D + A_{\tilde{D}},
\end{align*}
where,
\begin{align*}
    A_D &= \begin{pmatrix}
    A_{11} &  &  & \\
     & A_{22} &  & \\
     &  & \ddots & \\
     &  &  & A_{rr}\\
    \end{pmatrix},\\
    A_{\tilde{D}} &= \begin{pmatrix}
    0 & A_{12} & \cdots & A_{1r}\\
    A_{21} & 0 & \cdots & A_{2r}\\
    \vdots & \ddots & \ddots & \vdots\\
    A_{r1} & A_{r2} & \cdots & 0\\
    \end{pmatrix}.\\
\end{align*}

$A_D$ is called a \textit{pinching} of $A$. The simplest case is that of $A_{ij}$ being 1 dimensional that we will be using here.

Any positive semi-definite matrix $A$ can be written as $A= \sum_{i=1}^n \lambda_i P_{\lambda_i}$ where $\lambda_i$ are $n$ distinct eigen values  of $A$. $P_{\lambda_{i}}$ are orthogonal projectors such that $\sum_{i=1}^n P_{\lambda_i}=\mathbb{I}$  and hence $P_{\lambda_i}^2= P_{\lambda_i}$. The \textit{spectral pinching map} of $A$ is then given by
$$
\mathcal{P}_A: X \mapsto  \sum_{\lambda} P_\lambda XP_\lambda.
$$
The pinching map in turn has the following properties:
\begin{enumerate}[label=(\roman*)]
\item $\mathcal{P}_A[X] A = A\mathcal{P}_A[X]$.
\item $\mathrm{tr}[ \mathcal{P}_A[X]A] = \mathrm{tr} [XA]$. 
\item $\mathcal{P}_A[X] \geq  \frac{1}{n} X$.
\end{enumerate}

{\lem{} {$\mathcal{P}_A\left[X\right] A = A\mathcal{P}_A[X]$ }}.

\textit{Proof}:

Given that $P_{\lambda_i}^2 = P_{\lambda_i}$ and $P_{\lambda_i} \perp P_{\lambda_j}$,

\begin{align*}
    \mathcal{P}_A\left[X\right] A &= \sum_{i=1}^nP_{\lambda_i}XP_{\lambda_i} \sum_{i=1}^n\lambda_i P_{\lambda_i}\\
    &= \sum_{i=1}^n\sum_{j=1}^n\lambda_j P_{\lambda_i}XP_{\lambda_i} P_{\lambda_j}\\
    &= \sum_{i=1}^n\lambda_i P_{\lambda_i}XP_{\lambda_i} P_{\lambda_i}\\
    &= \sum_{i=1}^n\lambda_i P_{\lambda_i}XP_{\lambda_i}\\
    &= \sum_{i=1}^n\lambda_i P_{\lambda_i}P_{\lambda_i}XP_{\lambda_i}\\
    &= \sum_{i=1}^n\sum_{j=1}^n\lambda_j P_{\lambda_j}P_{\lambda_i}XP_{\lambda_i}\\
    &= \sum_{j=1}^n\lambda_j P_{\lambda_j} \sum_{i=1}^n P_{\lambda_i}XP_{\lambda_i}\\
    &= A\mathcal{P}_A[X].
\end{align*}
\qed

{\lem{} {$\mathrm{tr}\left[ \mathcal{P}_A[X]A\right] = \mathrm{tr} \left[XA\right]$ }}.

\textit{Proof}:
\begin{align*}
\mathrm{tr}\left[ \mathcal{P}_A\left[X\right]A \right] &= \mathrm{tr} \left[\sum_{i=1}^nP_{\lambda_i}XP_{\lambda_i} \sum_{i=1}^n\lambda_i P_{\lambda_i}\right]\\
 &= \mathrm{tr}\left[\sum_{i=1}^n\lambda_i P_{\lambda_i}XP_{\lambda_i}\right].\\
\mathrm{tr} \left[P_{\lambda_i}XP_{\lambda_i} \right] &= \sum_{j}(P_{\lambda_i}X)_{jj}\\
\implies \mathrm{tr}\left[ \mathcal{P}_A\left[X\right]A\right] &=  \sum_{i=1}^n\lambda_i\sum_{j}(P_{\lambda_i}X)_{jj} \\
&=\mathrm{tr}\left[AX\right].
\end{align*}
\qed

{\lem{} $\mathcal{P}_A\left[X\right] \geq  \frac{1}{n} X$ }.

\textit{Proof}:

Proving this part is probably the trickiest among the four lemmata here, but is the entire key behind deducing the final Golden-Thompson inequality. Consider a unitary matrix $U_y$ defined as $U_y = \sum_{u=1}^n e^{i2\pi yu/n}P_{\lambda_u}$. It is easy to verify that $U_yU_y^* =\mathbb{I}$ as $\sum_{i=1}^n P_{\lambda_i}=\mathbb{I}$. Also $U_y \succcurlyeq 0$ and $U_n=\mathbb{I}$

\begin{align*}
\sum_{y=1}^n U_yXU_y^* &= \sum_{y=1}^n\sum_{s=1}^n\sum_{t=1}^n e^{i2\pi ys/n}P_{\lambda_s}XP_{\lambda_t}e^{-i2\pi yt/n}\\
&= \sum_{s=1}^n\sum_{t=1}^n P_{\lambda_s}XP_{\lambda_t}\sum_{y=1}^ne^{i2\pi y(s-t)/n}\\
&= \sum_{s=1}^n\sum_{t=1}^n P_{\lambda_s}XP_{\lambda_t} (n\mathbf{1}_{\{s=t\}})\\
&= n\sum_{s=1}^nP_{\lambda_s}XP_{\lambda_s}\\
\implies \mathcal{P}_A[X] &= \sum_{s=1}^n P_{\lambda_s}XP_{\lambda_s} =  \frac{1}{n}\sum_{y=1}^nU_yXU_y^*\\
&\geq \frac{1}{n} X.
\end{align*}
\qed

Now, once Lemma 3 is proved, the rest of the steps for proving the Golden-Thompson inequality are straightforward. 
For a semi-positive definite $d \times d$ matrix $A$, we have the following lemma.

{\lem{}$|\mathrm{spec}(A^{\bigotimes m})| \leq O(\mathrm{poly}(m))$}

\textit{Proof}: 

The number of eigen values for $A^{\bigotimes m}$ is bounded by the number of the number of possible possible combinations of a sequence of $d$ symbols (the maximum possible distinct eigen values of $A$) of length $m$ which is given by $\binom{m+d-1}{d-1} \leq \frac{(m+d-1)^{d-1}}{(d-1)\,!} = O(\mathrm{poly}(m))$.

\clearpage

\section{Proof of Golden-Thompson inequality}

Given positive definite matrices $A, B$ and using the facts that $\exp{}$ and $\mathrm{tr}[\exp{}]$ are operator monotone:

\begin{align*}
    \log {\mathrm{tr}[\exp{\log A + \log B}]} &= \frac{1}{m} \log {\mathrm{tr}[\exp{\log A^{\bigotimes m} + \log B^{\bigotimes m}}]}\\
    &\leq  \frac{1}{m} \log {\mathrm{tr}[\exp{\log \mathcal{P}_{B^{\bigotimes m}}[A^{\bigotimes m}] |\mathrm{spec}(A^{\bigotimes m})| + \log B^{\bigotimes m}}]}  &&  \text{Using Lemma 3  }\\
    &=  \frac{1}{m} \log {\mathrm{tr}[\exp{\log \mathcal{P}_{B^{\bigotimes m}}[A^{\bigotimes m}] +  \log |\mathrm{spec}(A^{\bigotimes m})| + \log B^{\bigotimes m}}]}\\
    &=  \frac{1}{m} \log {\mathrm{tr}[\exp{\log \mathcal{P}_{B^{\bigotimes m}}[A^{\bigotimes m}] + \log B^{\bigotimes m}}]} +  \frac{1}{m}|\mathrm{spec}(A^{\bigotimes m})|\\
    &\leq  \frac{1}{m} \log {\mathrm{tr}[\exp{\log \mathcal{P}_{B^{\bigotimes m}}[A^{\bigotimes m}] + \log B^{\bigotimes m}}]} +  \frac{O(\mathrm{poly}(m))}{m}  && \text{Using Lemma 4} \\
    &= \frac{1}{m} \log {\mathrm{tr}[\exp{\log \mathcal{P}_{B^{\bigotimes m}}[A^{\bigotimes m}]B^{\bigotimes m}}]} +  \frac{O(\mathrm{poly}(m))}{m}  && \text{Using Lemma 1} \\
    &= \frac{1}{m} \log {\mathrm{tr}[\mathcal{P}_{B^{\bigotimes m}}[A^{\bigotimes m}]B^{\bigotimes m}]} +  \frac{O(\mathrm{poly}(m))}{m} \\
    &= \frac{1}{m} \log {\mathrm{tr}[A^{\bigotimes m}B^{\bigotimes m}]} +  \frac{O(\mathrm{poly}(m))}{m} && \text{Using Lemma 2}\\
    &= \log {\mathrm{tr}[AB]} +  \frac{O(\mathrm{poly}(m))}{m}\\
    \implies \log {\mathrm{tr}[\exp{\log A + \log B}]} &\leq \log {\mathrm{tr}(AB)}\\
    \implies  {\mathrm{tr}[\exp{\log A + \log B}]} &\leq {\mathrm{tr}(AB)}\\
    \implies  {\mathrm{tr}[\exp{ A + B}]} &\leq {\mathrm{tr}(\exp{A}\exp{B})}.\\
\end{align*}
\qed
\bibliographystyle{ieeetr}
\bibliography{references}
\end{document}

%% file: preamble.tex
\usepackage{amsmath}

\newcommand{\BIT}{\begin{itemize}}
\newcommand{\EIT}{\end{itemize}}
\newcommand{\BNUM}{\begin{enumerate}}
\newcommand{\ENUM}{\end{enumerate}}

\renewcommand{\exp}[1]{\operatorname{exp}\left(#1\right)} 

\newcommand{\interior}[1]{%
  {\kern0pt#1}^{\mathrm{o}}%
}

\newtheorem{lem}{Lemma}

\theoremstyle{remark}

\newtheorem*{rem*}{Remark}

\theoremstyle{definition}

\DeclareSIUnit\dbm{\decibel{}m}
\DeclareSIUnit\db{\decibel}

\newcolumntype{L}[1]{>{\raggedright\let\newline\\\arraybackslash\hspace{0pt}}m{#1}}
\newcolumntype{C}[1]{>{\centering\let\newline\\\arraybackslash\hspace{0pt}}m{#1}}
\newcolumntype{R}[1]{>{\raggedleft\let\newline\\\arraybackslash\hspace{0pt}}m{#1}}